\documentclass[11pt, a4 paper, reqno]{article}

\usepackage{multicol} 
\usepackage{makeidx}  
\usepackage{amsmath,amscd,amsfonts}
\usepackage{amssymb, fancyheadings,latexsym}

\newcommand{\Hm}{{\rm{Hom}}}
\newcommand{\F}{F_{\infty}}

\newcommand{\G}{G_{\infty}}

\newcommand{\kr}{{\rm{Ker}}}

\newcommand{\cok}{{\rm{Coker}}}

\newcommand{\Gal}{{\rm{Gal}}}

\newcommand{\cd}{{\rm{cd}}}

\newcommand{\q}{{\mathbb Q}}

\newcommand{\z}{\mathbb Z}
\newcommand{\f}{{\mathbb F}_p}
\newcommand{\lle}{\leqslant}
\newcommand{\gge}{\geqslant}
\newcommand{\fs}{F^S}
\newcommand{\rk}{{\rm{rank}}}
\newcommand{\C}{\mathcal{C}}
\newcommand{\D}{\mathcal{D}}
\newcommand{\E}{E_{p^{\infty }}}
\newcommand{\Sel}{{\cal S}}
\newcommand{\CC}{{\cal C}}

\newcommand{\ind}{{\rm{ind}}}
\newcommand{\res}{{\rm{res}}}
\newcommand{\hmrank}{{\rm{hmrank}}}
\newcommand{\ord}{{\rm{ord}}}

\newcounter{nonumber}

\newtheorem{thm}{Theorem}[section]
\newtheorem{lem}[thm]{Lemma}
\newtheorem{cor}[thm]{Corollary}
\newtheorem{prop}[thm]{Proposition}
\newtheorem{conj}[thm]{Conjecture}
\newtheorem{props}[nonumber]{Properties of {\it{\bfseries{p}}}-adic, Lie 
groups}
\newtheorem{wprops}[nonumber]{Properties}
\newcommand{\bigskp}{\par\addvspace{\bigskipamount}}
\newenvironment{proof}{\bigskp\noindent{\bf Proof }}{\bigskp}
\newenvironment{ass}{\bigskp\noindent{\bf Assumption }}{\bigskp}
\newenvironment{defn}{\bigskp\noindent{\bf Definition }}{\bigskp}
\newenvironment{ex}{\bigskip\noindent{\bf Example }}{\bigskip}
\newenvironment{Remark}{\bigskip\noindent{\bf Remark }}{\bigskip}
\newenvironment{Remarks}{\bigskip\noindent{\bf Remarks }}{\bigskip}

\newenvironment{question}{\bigskip\noindent{\bf Question }}{\bigskip}
\newenvironment{notn}{\bigskip\noindent{\bf Notation }}{\bigskip}

\begin{document}

\title{\vspace*{-2cm}Euler Characteristics as Invariants of Iwasawa Modules\\ \vspace*{-.3cm}}

\author{Susan~Howson{\setcounter{footnote}{6}\footnote{School of Mathematical Sciences, University of Nottingham, University Park, Nottingham, {\mbox{NG7 2RD}}, U.K.  {\ttfamily{\hspace*{.5em}sh@maths.nott.ac.uk}}  }} }
\date{}
\maketitle

\vspace*{-.8in}

\section*{}


\vspace*{-1cm}\section*{Introduction}
Let $G$ be a pro-$p$, $p$-adic, Lie group with no element of order $p$
and let $\Lambda (G)$ denote the Iwasawa algebra of $G$,  
defined in the usual way by
\begin{equation}\label{1}
\Lambda (G) = \underset{\longleftarrow}{\lim}\, \z _p [G/H],
\end{equation}
where $H$ runs over the open, normal subgroups of $G$, and the inverse 
limit is taken with respect to the canonical projection maps. 
Then many situations arise in which one is interested in determining 
information about the structure of modules which are finitely generated over $\Lambda (G)$, such as the 
rank (defined in (\ref{4}) below).
In this paper we describe a number of invariants associated to
a finitely generated $\Lambda (G)$ module, $M$, and calculated via an Euler
characteristic formula.
For example, we give a simple formula for the rank of $M$ in terms of an Euler Characteristics formula.  
These ideas are well known to Algebraists, see for example the book by K.~Brown, \cite{Brown}, chapter IX in particular, but do not appear to have been exploited to their full potential in Iwasawa Theory yet. This formula gives the natural 
generalisation 
of the strongest form of 
Nakayama's lemma (for the case of $G$ isomorphic to $\z _p$) to other
pro-$p$, $p$-adic, Lie groups. Thus the first subsection of 
this paper should be seen as a continuation 
of the earlier note \cite{BaHo97} which explained situations where 
that result can fail to generalise directly.

We then apply these ideas to modules finitely generated over the $\f $--linear completed group algebra, and consider an invariant of Iwasawa modules which gives the classical Iwasawa $\mu $-invariant in the case $G \cong \z _p $, where 
the idea of expressing this invariant in terms of an Euler characteristic is well known. 

If one instead starts with the Euler characteristic formula as the definition 
of a 'homological 
$\Lambda (G)$--rank' then we can extend it to some situations where the 
na\"{\i}ve definition of $\Lambda (G)$--rank is not 
appropriate. For example, we can consider removing 
the restriction that $G$ is a pro-$p$ group.
In the second section we give some discussion of the information this can 
tell us about a $\Lambda (G)$--module. In particular, in the classical case of 
$G \cong \z _p ^{\times}$ this relates to the decomposition of a finitely
generated $\Lambda (G)$--module into its eigenspaces for the action
of the subgroup ${\z / p\z }^{\times }$ of $\z _p ^{\times}$, 
a technique which is
not available to us for an arbitrary $p$-adic, Lie group. 

We also consider Euler characteristics of modules 
which are not actually finitely generated. In this way we can shed some 
light on a question raised in \cite{CoHo99} concerning a generalisation 
of Iwasawa's classical $\lambda $-invariant, which has proved 
important in the theory for $\z _p $--extensions of fields.  We consider in detail this invariant and the generalisation of the $\mu $-invariant for the dual of the Selmer group of an elliptic curve over a certain naturally defined $p$-adic, Lie extension, and analyse the behaviour of this second invariant under isogeny subject to some restrictions.

I would like to thank John Coates for many interesting comments and 
suggestions, R. Sujatha for pointing out that Proposition \ref{o}
can be deduced easily as a consequence of Shapiro's Lemma and the 
referee for many improvements to the exposition.

\section{Ranks of Iwasawa Modules}

\subsection{$ {\mathbb{\Lambda }}$--modules}\label{1.1}

Recall that $G$ denotes a $p$-adic, Lie group and $\Lambda (G)$ denotes the 
Iwasawa algebra of $G$, defined in (\ref{1}).
Let $M$ be any finitely generated $\Lambda (G)$--module. 
Then we are interested in the $\Lambda (G)$--rank of $M$. For the 
most natural definition of this
we need $\Lambda (G)$ to contain no non-trivial zero divisors. 
It is known 
that this is the case if $G$ contains no non-trivial element of finite order
(see \cite{Ne88})
and we assume this for the remainder of this section. In particular, this 
means $G$ is actually pro-$p$.
Then $\Lambda (G)$ is both right and left Noetherian and thus admits 
a skew field of fractions which we denote by $K(G)$ (see \cite{GoWa}, 
chapter 9.) 
\begin{defn} 
Under the assumption that $G$ contains no non-trivial element of finite order,
the $\Lambda (G)${\em{--rank}} of a finitely generated 
$\Lambda (G)$--module,
$M$, is defined by
\begin{equation}\label{4}
\rk _{\Lambda (G)}(M) = \dim _{K(G)} ( K(G)
\otimes _{\Lambda (G)} M).
\end{equation}
When this number is zero then we say $M$ is $\Lambda (G)$--{\em{torsion}}. If 
this is the case, then  
every element of $M$ has a non-trivial annihilator in $\Lambda (G)$, though it should be emphasised that in the non-Abelian setting this is not sufficient to ensure the existence of a non-trivial element in ann$_{\Lambda (G)} (M)$. 
\end{defn}
We will replace this later (equation (\ref{28})) 
with a less intuitive, homological definition of 
$\Lambda (G)$--rank which turns out to be better suited for generalising. In 
particular, it provides further insight into 
a question raised in \cite{CoHo99} 
which concerns modules that are not necessarily finitely generated.
Theorem \ref{a} of this section then states that 
the two definitions coincide for finitely generated 
$\Lambda (G)$--modules.

We start by recalling the natural duality between cohomology and homology.
If $M$ is a $\Lambda (G)$--module then we denote the Pontrjagin dual of $M$ by
\begin{equation}\label{22}
\widehat{M} = \Hm _{\z _p }(M,\q _p /\z _p ),
\end{equation}
which is naturally also a $\Lambda (G)$--module. Let $\C $ denote the 
category of finitely generated $\Lambda (G)$--modules, $\D $ the category of 
discrete, $p$--primary $\Lambda (G)$--modules (with continuous 
$\Lambda (G)$ action) which are cofinitely generated. We recall here
that a $\Lambda (G)$--module, $D$, is said to be cofinitely generated if 
$\widehat{D} $ is finitely generated. Then $\, \widehat{ } \, $ defines 
functors $\D \longrightarrow \C$ and $\C 
\longrightarrow \D$ 
(which for convenience we also denote by  $\; \widehat{ }$ ) cf. 2.3 of \cite{Br66}, which discusses this in greater generality. 
Taking the $G$--invariants, $D^G$, of a $\Lambda (G)$--module 
$D$ in $\D $ gives a left exact functor from $\D $ to the category of 
cofinitely generated $\z _p $--modules. Similarly, taking the 
$G$--coinvariants (the maximal quotient of $M$ on which $G$ acts trivially), 
$M_G$, of a $\Lambda (G)$--module $M$ in $\C $ gives a right exact functor 
from $\C $ to the category of finitely generated $\z _p $--modules. Since 
there is a canonical isomorphism
\begin{equation}\label{23}
\widehat{M_G} \cong \widehat{M} ^G
\end{equation}
for any $M$ in $\C $, this induces a canonical isomorphism between the 
corresponding derived  functors,
\begin{equation}\label{24}
H^i(G,\widehat{M}) \cong \widehat{H_i(G,M)},
\end{equation}
for any finitely generated $\Lambda (G)$--module, $M$. For convenience 
we prove our formula for modules in
the category $\C $, and 
use homology. Many applications in Number Theory, however, are 
phrased in terms of the $G$--cohomology of modules in $\D $ 
and thus we state our theorem in this language also.

Our first result is the following:
\begin{thm}\label{a}
Assume $G$ is a pro-$p$, $p$-adic, Lie group which 
contains no element of order $p$. Let
$M$ be a finitely generated $\Lambda (G)$--module.
Then the $\Lambda (G)$--rank of $M$ is given by the following 
'Euler Characteristic' formula:
\begin{equation}\label{2}
\rk _{\Lambda (G)}(M) = \sum_{i \gge 0} (-1)^{i} \,
\rk _{\z _p} \big( H_i(G,M) \big).
\end{equation}
Equivalently, for any $D$ in the category $\D $ of cofinitely generated, 
discrete $\Lambda (G)$--modules,
\begin{equation}\label{25}
\rk _{\Lambda (G)}(\widehat{D}) = \sum_{i \gge 0} (-1)^i \, 
{\rm{corank}} _{\z _p } \big( H^i (G,D) \big).
\end{equation}
In particular, both these sums are well-defined and finite.
\end{thm}
We observe that equality (\ref{25}) follows  
immediately from equality (\ref{2})
because of the duality in (\ref{24}). 
\begin{defn}
Let $G$ be a pro-$p$, $p$-adic, Lie group which contains no element of
order $p$. For any finitely generated $\Lambda (G)$-module, $M$, we 
define $\chi _* (G,M)$, respectively $\chi ^* (G,M)$, to be the integer
which is equal to the sum which occurs on the right hand side of equality 
(\ref{2}), respectively (\ref{25}).
\end{defn}

The following corollaries are immediate.
\begin{cor}\label{d}
Under the hypotheses of Theorem \ref{a},
if the $H_i(G,M)$ are finite for all $i \gge 0$ then $M$ is 
$\Lambda (G)$--torsion.
\end{cor}
In fact, since the $\Lambda (G)$--rank of any $M$ in $\C $ cannot be 
negative, it is sufficient to know the finiteness of $H_i (G,M)$ for $i$ 
even, only. As remarked above, it is interesting to compare this 
with the stronger results previously known in the classical case of 
$G \cong \z _p $, usually denoted by $\Gamma $. In this situation the 
following result is well known.
\begin{lem}\label{b}
Suppose $\Gamma $ is isomorphic to $\z _p$, and $M$ is a compact 
$\Lambda (\Gamma )$--module. If $M_{\Gamma }$ is a finite group then 
$M$ is $\Lambda (\Gamma )$--torsion.
\end{lem}
This extremely useful condition for a finitely generated 
$\Lambda (\Gamma )$--module to be torsion was first pointed out by Iwasawa 
and follows easily from the general structure theorem 
for finitely generated $\Lambda (\Gamma )$--modules. 
Note that it is immediate from this structure theorem 
that if $M_{\Gamma }$ is finite then so 
also is $M^{\Gamma }$. However, this now follows directly from 
Theorem \ref{a} above, without the crutch of the structure theorem.
Indeed, Lemma \ref{b} follows directly from Corollary \ref{d}, and the comment
which follows it, together with the equality 
$H_0 (\Gamma ,M) = M_{\Gamma }$ and and the fact that the cohomological 
dimension at $p$ of $\z _p $ is equal to 1. Then, since  
$H_1(\Gamma ,M) = M^{\Gamma }$, equality (\ref{2}) implies that 
$M^{\Gamma }$ is finite.
In \cite{BaHo97} it was shown that Lemma \ref{b} generalises 
to uniform, soluble, $p$-adic, Lie groups. However, as was explained 
there, Lemma \ref{b} fails for arbitrary insoluble groups, even under 
the assumption of being uniform. Thus Theorem \ref{a} seems 
to be the natural generalisation for arbitrary torsion free, pro-$p$, 
$p$-adic, Lie groups. See \cite{DdSMS} for a further discussion of the 
technical 
condition 'uniform'. We remark that $\z _p $ and also any 
congruence kernel of GL$_n (\z _p )$ are uniform groups.
As we will need the concept again later we define it here.
\begin{defn}
A profinite group is {\em{uniform}} if it is topologically finitely generated on $d$ generators, and there exists a filtration by subgroups
\begin{equation}
G = G_0 \supset G_1 \supset \cdots 
\end{equation}
such that each $G_i $ is normal in $G_{i+1}$, and $G_i / G_{i+1} \cong (\z / p\z ) ^d$. In particular, a uniform $p$-adic analytic group is always pro-$p$.
\end{defn}
\begin{Remark}
A fundamental fact is that a profinite group is $p$-adic analytic if and only if it contains an open, normal subgroup which is uniform \cite{DdSMS}.
\end{Remark}

\begin{cor}\label{e}
Assume the hypotheses of Theorem \ref{a}. 
If $M$ is finitely generated as a $\z _p $--module then
\begin{equation}\label{10}
\sum_{i \gge 0} (-1)^i \rk _{\z _p } \big( H_i (G,M) \big) = 0.
\end{equation}
\end{cor}
The only comment necessary here is to note that since we are assuming $G$ 
contains no element of order $p$, if it is not trivial (in which 
case Corollary \ref{e} is vacuous) then 
it must be infinite. Thus $\Lambda (G)$ is not finitely 
generated as a $\z _p$--module and so if $M$ is finitely generated as a 
$\z _p $--module then it must be $\Lambda (G)$--torsion.

The main application of these results in Number Theory is in Iwasawa Theory, 
where one is interested in studying the behaviour of arithmetic objects 
over infinite towers of fields. The typical situation is the following. Let 
$F$ be a finite extension of $\q $ and $\F $ a $p$-adic, Lie extension of 
$F$, that is a Galois field extension with $G=\Gal (\F /F)$ isomorphic to a 
$p$-adic, Lie group. Then certain modules in $\D $ are of great arithmetic 
interest. For example the cohomology groups
\begin{equation}\label{26}
H^i(F^S/\F ,A),
\end{equation}
where $F^S$ is the maximal extension of $F$ unramified outside a chosen 
finite set of primes of $F$, including all the primes which divide $p$
or which ramify in the 
extension $\F /F$ (assumed to be a finite set), 
and $A$ is isomorphic to ${({\q _p }/{\z _p})}^r $ as a 
$\z _p $--module, endowed with
a continuous action of $\Gal (F^S/F)$.

The application of Theorem \ref{a} to studying these modules, among other 
examples, will be described in \cite{Ho1}. Unfortunately, 
it is often hard to calculate higher homology and cohomology 
groups ($i \gge 2$) and so Theorem \ref{a} provides limited help, for example, 
with calculating the $\Lambda (G)$--rank of the $p^{\infty }$--Selmer 
group of an elliptic curve, as discussed in \cite{CoHo99}. We will also discuss
this situation in more detail later.

It is necessary to recall some general definitions and
properties of $p$-adic, Lie groups.

\begin{defn}
Recall that the $p$--{\em{cohomological dimension}} of $G$, denoted throughout by $\cd _p (G)$, 
is defined as the minimum number such that $H^i(G,D)$ vanishes 
for all discrete, 
$p$--primary, $\Lambda (G)$--modules, $D$, whenever $i > \cd _p (G)$.
It follows from (\ref{24}) that for all finitely generated $\Lambda (G)$--modules, $M$, one also has 
$H_i (G,M)=0$ for all  
$i > \cd _p (G)$.

If $M$ is in $\C $ then the 
$\Lambda (G)$--{\em{homological dimension}} of $M$ is the 
minimum integer, $n$, 
such that $M$ has a resolution of length $n$ by projective
$\Lambda (G)$--modules. (The $\Lambda (G)$--homological dimension is said to 
be infinite if no finite resolution exists.) We define the 
{\em{global dimension}}, $gldim(\Lambda (G))$, of $\Lambda (G)$ to be 
the supremum of the homological dimension of $M$, as $M$ ranges over $\C $.
\end{defn}

\begin{props}\label{c}
\begin{enumerate}
\item Any $p$-adic, Lie group containing no element of order $p$ (but 
which is not
necessarily pro-$p$)
has finite $p$--cohomological dimension, 
equal to its dimension as a $p$-adic Lie group, \cite{La65} and \cite{Se71-2}.
\item It is shown in \cite{Br66} 
that under the hypotheses of i) $\Lambda (G)$ has finite global dimension, 
given by
\begin{equation}\label{3}
gldim(\Lambda (G)) = 1 + \cd _p (G).
\end{equation}
\item For $G$ pro-$p$, but not necessarily torsion free, 
$\Lambda (G)$ is a local ring, with unique maximal ideal given by the 
kernel of the augmentation map, $\Lambda (G) \rightarrow \z / {p\z }$, 
\cite{La65}.
\end{enumerate}
\end{props}

\begin{proof}of Theorem \ref{a}.\\
We fix $G$ throughout the proof of Theorem \ref{a} 
and write $\Lambda $ for $\Lambda (G)$.
We also fix a finitely generated 
$\Lambda $--module, $M$.
The finiteness of the global dimension of $\Lambda $ means that
$M$ admits a projective resolution of finite length (at most equal to 
$gldim(\Lambda )$, by definition.) Because $\Lambda $ is a 
local ring, the only projective $\Lambda $--modules are the free modules. 
Thus $M$ admits a resolution of finite length, $d$ say,
by free modules. Fix such a resolution:
\begin{equation}\label{5}
0 \rightarrow \Lambda ^{n_d} \rightarrow \cdots \rightarrow \Lambda ^{n_1} 
\rightarrow \Lambda ^{n_0} \rightarrow M \rightarrow 0.
\end{equation}
It follows that the groups $H_i(G,M)$ can be computed as the homology of a complex of finite length of the form
\begin{equation}
\z _p^{n_d} \longrightarrow \cdots \longrightarrow 
\z _p^{n_1} \longrightarrow \z _p ^{n_0 }
\end{equation}
immediately implying that each $H_i(G,M)$ is finitely generated over $\z _p $ and that 
\begin{equation}
\chi _* (G,M) = \underset{i \gge 0 } \sum (-1)^i n_i .
\end{equation}

On the other hand, since $K(G)$ is a flat $\Lambda $--module, 
the alternating sum of $\Lambda $--ranks adds to zero along exact 
sequences. Thus it follows from (\ref{5}) that 
\begin{equation}\label{17}
\rk _{\Lambda }(M) = n_0 - n_1 + n_2 \cdots + (-1)^d n_d ,
\end{equation}

completing the proof of Theorem \ref{a}.{\\ \nopagebreak
\hspace*{\fill}$\Box$}
\end{proof}

The following simple Corollary will be needed later.
\begin{cor}\label{qq}
If $G_0 $ is an open normal subgroup in $G$ then the rank of a finitely generated $\Lambda (G)$--module over $\Lambda (G_0 )$ equals $ | G:G_0 | $ times its rank over $\Lambda (G ) $. In particular, a $\Lambda (G)$--module is $\Lambda (G)$--torsion if and only if it is $\Lambda (G_0 )$--torsion.
\end{cor}
\begin{proof}
This is immediate from the fact that
$\Lambda (G)$ is free over $\Lambda (G_0 )$, of rank $ | G : G_0 | $, and so a free resolution of $M$ as a $\Lambda (G)$--module also gives a free resolution of $M$ as a $\Lambda (G_0 )$--module, with alternating sum of $\Lambda (G_0 )$--ranks equal to $| G:G_0 |$ times the alternating sum of $\Lambda (G)$--ranks. 
{\\ \nopagebreak
\hspace*{\fill}$\Box$}
\end{proof}

\begin{ex}
As a simple illustration of the application of Theorem \ref{a} we consider 
the smallest concrete example discussed in \cite{BaHo97} of a $p$-adic, Lie
group for which the 
assertion of Lemma \ref{b} can fail.

Assume $p \gge 5$ and let
\begin{equation}\label{50}
H = \big\{ A \in {\rm{SL}}_2 (\z _p ) \mid A \equiv I ,\; {\rm{mod}}\;p \big\}.
\end{equation}
Then $H$ is not soluble (that is, its derived series does not terminate 
after a finite number of terms) and it can be shown that 
\begin{equation}\label{51}
\overline{[H,H]}=
\big\{ A \in {\rm{SL}}_2 (\z _p ) \mid A \equiv I ,\; {\rm{mod}}\; p^2 \big\},
\end{equation}
where $\overline{[H,H]}$ is the group generated 
{\em{topologically}} by $[H,H]$.
Now $H$ is a pro-$p$, $p$-adic, Lie group of dimension 3. Since $p \gge 5$ it 
contains no element of order $p$, and so Theorem \ref{a} applies.
Let $I(H)$ denote the augmentation ideal of $\Lambda (H)$, that is 
the kernel of the canonical augmentation map, $\varepsilon $, to $\z _p $
\begin{equation}\label{52}
0 \longrightarrow I(H) \longrightarrow \Lambda (H) \overset{\varepsilon }{\longrightarrow} 
\z _p \longrightarrow 0.
\end{equation}
Since the alternating sum of $\Lambda (H)$--ranks adds along exact sequences it is clear that
\begin{equation}\label{53}
\rk _{\Lambda (H)}(I(H)) =1.
\end{equation}
We now calculate the homology groups $H_i(H,I(H) )$, in order to illustrate 
that the non-zero contribution to the sum in the formula
(\ref{2}) can come entirely 
from the terms with $i \not= 0$. We start by taking homology of (\ref{52}).
This gives
\begin{equation}\label{54}
\begin{array}{rcll}
H_0(H,I(H)) &\cong &  I(H)/{(I(H))^2}\\
H_i(H,I(H))  &\cong  &  H_{i+1}(H,\z _p ),&\quad i \gge 1 \\
&  = & 0,  &\quad   i \gge 3,
\end{array}
\end{equation}
since $cd _p (H) =3$. It is easy to see that as an Abelian group
\begin{equation}\label{55}
I(H)/{((H))^2} \cong H/{\overline{[H,H]}},
\end{equation}
which is isomorphic to $(\z /{p \z})^3 $ by (\ref{50}) and (\ref{51}) 
and thus in particular is finite. 
To calculate $H_i(H,\z _p )$ for $i = 2,3$ we 
use the fact, shown in \cite{La65} Theorem 2.5.8, 
that any pro-$p$, $p$-adic, Lie group is 
a Poincar\'{e} group of dimension equal to its dimension as a 
$p$-adic manifold. The definition of a Poincar\'{e} group is somewhat 
technical and we simply recall the following basic property:

If $G$ is a Poincar\'{e} group of dimension $n$ then the 
cup product induces a perfect pairing
\begin{equation}\label{56}
H^i(G,A) \times H^{n-i} (G,\widetilde{A}) \longrightarrow \q _p /\z _p
\end{equation}
for $A$ any finite $G$--module endowed with the discrete topology. Here 
$\widetilde{A} $ denotes $\Hm _{cts} (A,I)$, endowed with the 
natural action of $G$, 
where $I$ is the dualising module of $G$. We shall not define 
dualising modules in general. In fact, $H$ is what is 
known as an {\it{orientable}} Poincar\'{e}
group, in other words its dualising module is simply 
$\q _p / \z _p $ with trivial action of $H$.
See \cite{La65} and \cite{Se1} for a more detailed 
discussion of the definition and properties of Poincar\'{e} groups. 

From the corollary to Proposition 2.2 of \cite{Ta76},
\begin{equation}\label{57}
H^i_{cts}(H,\z _p ) = \underset{\longleftarrow}{\lim} \,H^i (H,\z /p^n \z ),
\end{equation}
where the inverse limit is taken with respect to the canonical maps.
Here $H_{cts}^i(H,\z _p )$ denotes cohomology
calculated with continuous cochains and 
where $\z _p $ is endowed with the usual profinite topology. It
follows from this, the isomorphism in 
(\ref{24}) and the pairing (\ref{56}) above that
\begin{equation}\label{58}
H_i(H,\z _p) \cong H^{3-i}_{cts}(H,\z _p ).
\end{equation}
Since $H^1_{cts}(H,\z _p )$ is just the group of continuous group homomorphisms
of $H$ into $\z _p $, it follows from (\ref{58}) and the finiteness of 
$H^{ab} = H/\overline{[H,H]}$ that $H_2 (H,\z _p ) $ vanishes and so $H_1 (H,I(H))$ vanishes also.
Whereas putting together (\ref{54}) and (\ref{58}) for $i$ equal to 2, we see
\begin{equation}\label{60}
H_2(H,I(H)) \cong H_3(H,\z _p ) \cong H^0_{cts}(H,\z _p ) \cong \z _p.
\end{equation}
So the only non-zero contribution to the formula
(\ref{2}) is from the term $H_2(H,I(H))$, of $\z _p $--rank 
equal to 1.

\end{ex}

\subsection{${\mathbb{\Omega }}$--modules}\label{1.2}

We include here a short subsection concerning $\f $--linear, completed group algebras. This is interesting with regard to our attempts to understand more generally the possible structures of modules over Iwasawa algebras. Otmar Venjakob has independently considered the structure of modules over the $\f $--linear, completed group algebra in some detail in his Heidelberg PhD. thesis, \cite{Ven}.

We use the notation
\begin{equation}\label{200}
\Omega (G) =  \underset{\longleftarrow}{\lim}\, \f  [G/H],
\end{equation}
for the $\f $--linear, completed group algebra, and continue with the assumption that $G$ contains no non-trivial element of finite order. Clearly $\Omega (G) = \Lambda (G) /{p \Lambda (G) }$ and the following properties of $\Omega (G)$ hold, as for $\Lambda (G)$. 

\begin{wprops}
\begin{enumerate}

\item The global dimension of $\Omega (G)$ is finite, this time equal to the cohomological dimension of $G$ at $p$, since $\f $ is a field.
\item $\Omega (G)$ is a right and left Noetherian, local ring.
\item As remarked just prior to the statement of Corollary \ref{e}, $G$ contains an open normal subgroup, $G_0 $, which is uniform. The associated group algebra $\Omega (G_0 )$ has no non-trivial zero divisors (see the second edition of \cite{DdSMS}), thus admits a skew field of fractions, which we denote by $k (G_0 )$. 
\end{enumerate}
\end{wprops}
Unfortunately, the proof in \cite{Ne88} that $\Lambda (G)$ contains no non-trivial zero divisors when $G$ is torsion free fails to generalise to characteristic $p$, and it appears to be unknown whether $\Omega (G)$ is torsion free without the restriction that $G$ be uniform.

\begin{defn}
Similarly to the case for $\Lambda (G)$--modules, the $\Omega (G_0)${\em{--rank}} of a finitely generated $\Omega (G)$--module, $M$, is given by
$ \dim _{k(G_0)} ( k(G_0) \otimes _{\Omega (G_0 )}M) $.
Since $\Omega (G)$ is finitely generated over $\Omega (G_0 )$, $M$ is also finitely generated over $\Omega (G_0 )$. Thus we extend this definition and define the $\Omega (G)$--{\em{rank}} as 
\begin{equation}\label{201}
\rk _{\Omega (G)}(M) = \frac{\rk _{\Omega (G_0)} ( M )  }{ | G : G _0 | }.
\end{equation}
We will show below that this is integral and independent of the choice of $G_0 $.

Throughout the remainder we will use the notation
\begin{equation}\label{502}
\chi (G,M) =  \underset{i \gge 0}{\prod}\, \big( \# H_i (G,M) \big) ^{(-1)^i}
\end{equation}
for any $\Lambda (G)$--module, $M$, for which the terms in this product are finite and almost all equal to one.
\end{defn}
Then the arguments of section \ref{1.1} all follow through also for finitely generated $\Omega (G)$--modules. In particular, one has the following result.
\begin{prop}\label{z} 
Assume $G$ is a pro-$p$, $p$-adic, Lie group which contains no element of order $p$. Let $M$ be a finitely generated $\Omega (G)$--module. Then the $\Omega (G)$--rank of $M$ is given by the following, finite, sum
\begin{equation}\label{202}
\rk  _{\Omega (G)}(M) =  \sum_{i \gge 0} (-1)^{i} \,
\dim _{\f } \big( H_i(G,M) \big).
\end{equation}
By the definition above, this is equal to $\ord _p (\chi (G,M))$, where $\ord _p (\alpha )$
 denotes the maximum power of $p$ which divides $\alpha$. 
In particular, for any finitely generated $\Omega (G)$--module, $M$, it follows
that $\chi (G,M)$ is always well-defined and that the definition (\ref{201})
of the $\Lambda (G)$--rank is both independent of the choice of subgroup, $G_0 $, used to define it and integral.
\end{prop}
\begin{proof}
This follows by almost exactly the same arguments used above to prove 
Theorem \ref{a}. Since $\Omega (G)$ is a local ring and has finite 
global dimension, we may again take a finite, free resolution of $M$, 
this time as an $\Omega (G)$--module,
\begin{equation}\label{500}
0 \longrightarrow \Omega(G)^{n_d} \longrightarrow \cdots \longrightarrow \Omega (G) ^{n_0 } \longrightarrow M  \longrightarrow 0.
\end{equation}
The only extra point needed is that since $\Omega (G)$ is not a 
projective $\Lambda (G)$--module, we must first 
calculate $\chi (G,\Omega (G) ^n ) = \chi (G,\Omega (G) )^n $. 
Consider the sequence
\begin{equation}
0 \longrightarrow \Lambda (G) \overset{\times p }{\longrightarrow } \Lambda (G) \longrightarrow \Omega (G) \longrightarrow 0.
\end{equation}
Taking $G$--homology, and recalling $H_i (G,\Lambda (G) ) = 0 $ for $i \gge 1 $, we obtain
\begin{equation}
0 \longrightarrow H_1(G,\Omega (G)) \longrightarrow \z _p \overset{\times p }{\longrightarrow } \z _p \longrightarrow H_0 (G,\Omega (G)) \longrightarrow 0
\end{equation}
and $H_i (G,\Omega (G))$ vanishes for $i \gge 2 $. Since $H_1 (G,\Omega (G))$ is both annihilated by $p$ and contained in $\z _p $, it must be zero also, while $H_0 (G,\Omega (G)) $ is simply $\f $. Thus $\chi (G,\Omega (G)^n ) $ equals $n $ and similarly $\chi (G_0 , \Omega (G)^n ) $ equals $ | G : G_0 | n $.
It follows that we can calculate the $H_i(G_0,M)$ by working in the category of finitely generated $\Omega (G_0 )$--modules and using the projective resolution (\ref{500}).
The remainder of the proof of Proposition \ref{z} now follows exactly as before for Theorem \ref{a}, in particular $\Omega (G_0 )$ is again flat as an $\Omega (G_0)$--module and so the ${\Omega (G _0 )} $--rank of $M$ is  equal to $(n_0 - n_1 + \cdots \pm n_d )\mid G:G_0 \mid $, from which the proposition is clear.{\\ \nopagebreak
\hspace*{\fill}$\Box$}
\end{proof}

Corollaries \ref{d} and \ref{e} also have equivalent formulations for finitely generated $\Omega (G)$--modules. In fact the analogue of Corollary \ref{e} is already contained in the exercise at the end of \S I.4 of \cite{Se1} and is valid for more general groups, $G$. We will therefore omit reformulating statements or proofs.

Suppose now that $M$ is a finitely generated $\Lambda (G)$--module. Let $M (p)$ denote the subset of $M$ consisting of the elements of $M$ which are annihilated by some power of $p$ (clearly a $\Lambda (G)$--submodule as $p$ is in the centre of $\Lambda (G)$.) Since $\Lambda (G)$ is Noetherian, $M (p)$ is again finitely generated, and thus there exists some integer $r \gge 0$ such that $p^r $ annihilates $M(p)$. We propose the following invariant of $\Lambda (G)$--modules.
\begin{defn}
If $M$ is a finitely generated $\Lambda (G)$--module then we define
\begin{equation}\label{203}
\mu (M) = \sum_{i \gge 0}  \,\rk _{\Omega (G) } \Big( {p^i \big( M(p) \big)}/{p^{i+1}} \Big)
\end{equation}
\end{defn}
Taken in conjunction with the discussion preceding this definition, the fact
that each sub-quotient $\Big( {p^i \big( M(p) \big)}/{p^{i+1}} \Big)$ is a finitely generated $\Omega (G)$--module implies that this sum is finite. We have used the notation $\mu (M) $ since in the case $ G \cong \z _p $ this gives the classical $\mu $-invariant of Iwasawa Theory for $\Gamma $--modules.  This is a particularly convenient definition to work with.
\begin{cor}\label{dd}
Assume $G$ contains no non-trivial element of finite order. If $M$ is a finitely generated $\Lambda (G)$--module then
\begin{equation}\label{204}
\mu (M) = \ord _p \big( \chi (G,M(p)) \big) .
\end{equation}
In particular, $\chi (G, M(p))$ is always finite.
\end{cor}

Unfortunately, there exist finitely generated $\Lambda (G)$--modules,
$M$, for which $\chi (G,M)$ is finite and non-zero and yet $\mu (M)$ is equal to zero, and so we $\mu (M)$ cannot be definited solely in terms of $\chi (G,M)$. Indeed, as has been shown in \cite{CSW00}, there actually exist modules $T$ which are finitely generated and free over $\z _p $ and for which $\ord _p (\chi (G,T) ) $ is strictly positive, or even strictly negative.
These examples are obtained by twisting by roots of unity the Tate module of Tate curves defined over a finite extension of $\q _p $. See \cite{CSW00}, Corollary 5.3, for further details. In particular, $\chi (G,M)$ need not be integral if $M$ is not $p$-torsion, even for rather small $\Lambda (G)$--modules, and allowing $G$ to be uniform.
\begin{proof}
We consider the short exact sequences
\begin{equation}\label{205}
0 \longrightarrow p^{i+1} (M(p)) \longrightarrow p^i (M(p)) \longrightarrow
{p^i (M(p))}/{p^{i+1}} \longrightarrow 0
\end{equation}
obtained from the filtration
\begin{equation}\label{206}
M(p) \supset p(M(p)) \supset \cdots \supset p^r(M(p)) = 0
\end{equation}
for $r$ sufficiently large. The modules ${p^i (M(p))}/{p^{i+1}}$ are finitely generated $\Omega (G)$--modules for all $i$ and so, by Proposition \ref{z}
\begin{equation}\label{207}
\ord _p \chi \big( G,{p^i (M(p))}/{p^{i+1}} \big) = \rk _{\Omega (G)} \big( {p^i (M(p))}/{p^{i+1}} \big).
\end{equation}
In particular, these numbers are all defined. Since $p^r (M(p))$ vanishes, the Euler characteristic $\chi (G,p^{r-1}(M(p)))$ is finite by Proposition \ref{z}. It then follows inductively from the long exact sequence in cohomology that $\chi (G, p^i (M(p)))$ is finite for all $i$, equal to the product of $\chi \big( G, p^{i+1} (M(p)) \big)$ and $ \chi \big( G,{p^i (M(p))}/{p^{i+1}} \big)$. Multiplying these together, we obtain the formula (\ref{204}) of the 
Corollary.{\\ \nopagebreak
\hspace*{\fill}$\Box$}
\end{proof}

Although this invariant is clearly not additive along exact sequences (consider, for example, the sequence $0 \rightarrow p \Lambda (G) \rightarrow \Lambda (G) \rightarrow \Omega (G) \rightarrow 0 $) it does behave well when restricted to torsion $\Lambda (G) $--modules.

\begin{prop}\label{bb}
Assume $G$ contains no non-trivial element of finite order. 
In a short exact sequence of finitely generated $\Lambda (G)$--modules,
\begin{equation}
0 \longrightarrow A \longrightarrow B \longrightarrow C \longrightarrow 0,
\end{equation}
one has $\mu (B) \lle \mu (A) + \mu (C)$, with equality if $B$, and hence also $A$ and $C$, is $\Lambda (G)$--torsion.
\end{prop}
\begin{proof}
We start with the observation that if $A, B $ and $C$ are $p$--torsion then the lemma follows from Corollary \ref{dd}, since all Euler characteristics are defined and alternating products of Euler characteristics along exact sequences multiply together to give 1.
In general we have $0 \rightarrow A(p) \rightarrow B(p) \rightarrow C(p) $, which, together with this remark gives the inequality. However, exactness can fail on the right. Let $X$ denote the image of $B(p)$ in $C(p)$. It follows from the first remark that $\mu (B) = \mu (A) + \mu (X)$. We show that if $B$ is $\Lambda (G)$--torsion, then $\mu (C(p) / X) =0$ and thus $\mu (C) = \mu (X)$. 
The surjection $B \twoheadrightarrow C$ induces a surjection $B/B(p) \twoheadrightarrow C/X$. By assumption, $B$ is $\Lambda (G)$--torsion, and thus $B/B(p)$ is also $\Lambda (G )$--torsion.
By Corollary \ref{qq} it is also $\Lambda (G_0 )$--torsion for any choice of uniform, open, normal subgroup of $G$. Let $G_0 $ be any such subgroup. 
Then $B/B(p)$ is $\Lambda (G_0 )$--torsion and  $p$--torsion free. Since $p$ is in the centre of $\Lambda (G_0)$, it follows that for 
every element $c$ of $ C/X$ there exists an element, $\lambda $, of $ \Lambda (G_0)$ which is not contained in $p(\Lambda (G_0))$ and which annihilates $c$. In particular, this is true for all elements of $C(p) /X$.
\begin{lem}\label{cc}
Let $G_0 $ be a uniform, pro-$p$ group.
Suppose $Z = \Lambda (G_0) \alpha $ is a $p$--torsion $\Lambda (G_0)$--module, generated by a single element, $\alpha $, and such that there exists an element $\lambda $ of $  \Lambda (G_0) $, not contained in
$ p\Lambda (G_0) $, for which  $\lambda \alpha = 0$. Then $\mu (Z) = 0$.
\end{lem}
\begin{proof}
Since $p$ is in the centre of $\Lambda (G_0)$ we may consider the modules $p^i Z / p^{i+1} $, equal to $\Lambda (G_0)p^i \alpha / p^{i+1}$, separately for each $i$ and thus, without loss of generality, we may assume $Z$ is annihilated by $p$, that is, $Z$ is an $\Omega (G_0)$--module generated by a single element, $\alpha $. Furthermore, this generator $\alpha $ has a non trivial annihilator in $\Omega (G_0)$ given by the image of $\lambda $ under the canonical surjection of  $ \Lambda (G_0) $ onto $ \Omega (G_0)$. But then it is immediate this means  $Z$ has rank zero as an $\Omega (G_0)$--module, as required.{\\ \nopagebreak\hspace*{\fill}$\Box$}
\end{proof}
By the definition of $\Omega (G)$--ranks, for any finitely generated $\Omega (G)$--module, $N$, the $\mu $-invariant, $\mu (N)$,  calculated as an $\Omega (G)$--module vanishes if and only if it vanishes when considering $N$ as an 
$\Omega (G_0 )$--module. 
The module $C(p)/X $ is finitely generated over $\Lambda (G_0 )$, as $\Lambda (G_0)$ is Noetherian and $\Lambda (G)$ is finitely generated over $\Lambda (G_0 )$. Thus the proposition follows by induction on the number of generators of $C(p)/X$.{\\ \nopagebreak
\hspace*{\fill}$\Box$}
\end{proof}

We now finish this section with a result which compares the 
$\Omega (G)$--rank of $M/p$ with the $\Lambda (G)$--rank of $M$.
\begin{cor}
Assume $G$ contains no non-trivial element of finite order. 
Let $M$ be a finitely generated $\Lambda (G)$--module. Then 
\begin{equation}\label{300}
\rk _{\Omega (G)} (M/pM) = \rk _{\Omega (G)} (M[p]) + \rk _{\Lambda (G)} (M),
\end{equation}
where $M[p]$ denotes the submodule of $M$ consisting of the set of elements, $m$, contained in $M$ such that $pm =0$.
\end{cor}
\begin{proof}
We begin by considering a finitely generated $\Lambda (G)$--module, $N$, 
which is actually $p$--torsion. 
Since $N$ is $p$--torsion, it follows from Corollary \ref{dd} 
that $\chi (G,N)$,
$\chi (G,N[p])$ and $\chi (G,N/pN)$ are all finite, where $\chi (G, \_)$
is the Euler characteristic in the category of $\Lambda (G)$--modules as
defined by (\ref{502}). Upon consideration of the exact
sequence :
\begin{equation}\label{501}
0 \longrightarrow N[p] \longrightarrow N \longrightarrow N \longrightarrow N/pN \longrightarrow 0,
\end{equation}
and since the alternating product of Euler characteristics along an exact sequence multiplies together to make 1, this shows that 
\begin{equation}
\chi (G,N[p]) = \chi (G,N/pN ).
\end{equation}
However, both $N[p]$ and $N/pN $ actually have exponent $p$ and so, by Proposition \ref{z}, this shows that they have equal rank as $\Omega (G)$--modules.
Thus (\ref{300}) holds for 
finitely generated $p$--torsion modules.
In particular, it holds for $M(p)$. 

Because $M/M(p)$ is $p$--torsion free, the sequence 
\begin{equation}
0 \longrightarrow M(p)/pM(p) \longrightarrow M/pM \longrightarrow \frac{M/M(p)}{p \big( M/M(p) \big) } \longrightarrow 0
\end{equation}
is exact, and so
\begin{equation}
\rk _{\Omega (G)} (M/pM) = \rk _{\Omega (G)} \Big( \frac{M/M(p)}{p \big( M/M(p) \big) } \Big)  + \rk _{\Omega (G)} (M[p]).
\end{equation}
Since the $\Lambda (G)$--rank of $M/M(p)$ equals that of $M$, it only remains to show that if $N$ is a $p$--torsion free $\Lambda (G)$--module, then the $\Omega (G)$--rank of $N/pN$ equals the $\Lambda (G)$--rank of N. But this follows immediately from Theorem \ref{a}, Proposition \ref{z} and the long exact sequence in $G$-homology of 
\begin{equation}
0 \longrightarrow N \overset{\times p }{\longrightarrow } N \longrightarrow N/p \longrightarrow 0.
\end{equation}
Indeed, the multiplication by $p$ map on N induces multiplication by $p$ on the homology groups, $H_i(G,N)$. Since on a finite Abelian group the kernel and cokernel of multiplication by $p$ have the same order, we easily see that
\begin{equation}
\ord _p (\chi (G,N/p) )= \chi _* (G,N).
\end{equation}
By Theorem \ref{a}, the term on the right equals the $\Lambda (G)$--rank of $N$, while by Proposition \ref{z} the term on the left equals the $\Omega (G)$--rank of $N/pN$.{\\ \nopagebreak
\hspace*{\fill}$\Box$}
\end{proof}

\section{Homological $\mathbb{\Lambda }$--ranks}

We have proven that under the conditions of Theorem \ref{a} 
the right hand side of 
(\ref{2}) is always finite. However, there are other situations where $\chi _* (G,M)$ can be finite and the na\"{\i}ve definition of rank is no
longer suitable. Thus we propose the following definition.
\begin{defn}
Let $G$ be any $p$-adic, Lie group with no elements of order $p$, and 
$\Lambda (G)$ the Iwasawa algebra of $G$. Let $M$ be any compact 
$\Lambda (G)$--module. Then we say $M$ has {\em{finite homological rank}}
if the sum on the right of (\ref{2}), $\chi _* (G,M)$, has finitely many non-zero, but finite, terms and is thus well-defined. In which case we define
\begin{equation}\label{28}
{\rm{hmrank}}_{\Lambda (G)}(M) = \chi _* (G,M) . 
\end{equation}
\end{defn}
So Theorem \ref{a} states that, for $G$ a pro-$p$, $p$-adic, Lie group, containing no element of order $p$,
\begin{equation}\label{29}
{\rm{hmrank}}_{\Lambda (G)} (M) = \rk _{\Lambda (G)} (M),
\end{equation}
if $M$ is finitely generated.

It follows from the long exact sequence in $G$--homology that
\begin{lem}\label{g}
homological rank is additive along exact sequences. More precisely for any $p$-adic, Lie group, let
\begin{equation}\label{30}
0 \longrightarrow M \longrightarrow N \longrightarrow P \longrightarrow 0
\end{equation}
be an exact sequence of compact $\Lambda (G)$--modules. Then if any 
two have finite homological rank, so does the third, in which case
\begin{equation}\label{31}
{\rm{hmrank}}_{\Lambda (G)} (N) = {\rm{hmrank}}_{\Lambda (G)} (M) + 
{\rm{hmrank}}_{\Lambda (G)} (P).
\end{equation}
\end{lem}

We believe that this invariant will turn out to be very useful in a more general structure theory for $\Lambda (G)$--modules. In particular, as evidence we consider the following two situations:

\subsection{{\it{\bfseries{G}}} fails to be pro-{\it{\bfseries{p}}}}
We can now consider $p$-adic, Lie groups, $G$, which are not 
pro-$p$. This is necessary for the groups which turn up 'in nature', for 
example $\Gal (\q (\mu _{p^{\infty }}) /\q )$ or $\Gal (\q (\E )/\q )$,
where $\mu _{p^{\infty }} $ is the set of all $p^{\rm{th}}$--power roots 
of unity and $\q (\E )$ is the field of definition of the 
$p^{\rm{th}}$--power torsion points on an elliptic curve, $E$, 
defined over $\q $. Both of these field extensions contain the extension 
$\q (\mu _p )/\q $ 
which has degree $p-1$. The na\"{\i}ve definition of $\Lambda (G)$--rank 
does not generalise. Indeed, if $G$ is not a pro-$p$ group then it will 
contain non-trivial elements 
of finite order and so $\Lambda (G)$ will contain zero divisors other 
than the zero element.
Thus there is no skew field into which $\Lambda (G)$ can be embedded.
\begin{lem}\label{h}
If $M$ is a finitely generated $\Lambda (G)$--module, for $G$ a 
$p$-adic, Lie group containing no element of order $p$, 
then the homological rank of $M$ is always finite.
\end{lem}
\begin{proof} 
From the remark just prior to the statement of Corollary \ref{e}, we know
$G$ contains a pro-$p$, 
normal subgroup of finite index, Thus the finiteness of the $\z _p $--rank
of $H_i(G,M) $ follows immediately from the case of a pro-$p$, $p$-adic, Lie
group and the 
Hochschild-Serre spectral sequence for group cohomology.{\\ \nopagebreak
\hspace*{\fill}$\Box$}
\end{proof}

It is not clear how to interpret this homological dimension, but an 
indication is the following.

\begin{prop}\label{o}
Let $G$ be a $p$-adic, Lie group containing no element of order $p$.
Let $G_0 $ be an open subgroup of $G$ which is pro-$p$. 
Suppose $M$ is a finitely generated 
$\Lambda (G)$--module. Let $\res ^G _{G_0} M $ denote $M$ considered as a $\Lambda (G_0 )$--module by restricting the action of $\Lambda (G)$. For any finitely generated $\Lambda (G_0 )$--module, $N$, let $\ind ^G _{G_0 } N = N \otimes _{\Lambda (G_0 )} \Lambda (G)$ be the algebraically induced $\Lambda (G)$--module. Then 
\begin{equation}\label{210}
\hmrank _{\Lambda (G) } (\ind  ^G _{G_0 } ( \res ^G _{G_0}M)) = \rk _{\Lambda (G_0 )} (M).
\end{equation}
In particular, if $G_0 $ is normal in $G$, with $ \Delta = G/G_0 $ then
\begin{equation}\label{211}
\hmrank _{\Lambda (G) } (M \bigotimes _{\z _p } \z _p [ \Delta ] ) = \rk _{\Lambda (G_0 ) } (M)
\end{equation}
where $\z _p [\Delta ] $ is a $\Lambda (G) $--module via the canonical map $\Lambda (G) \rightarrow \z _p [\Delta ] $.
\end{prop} 
\begin{proof}
This is an immediate consequence of Shapiro's Lemma together with Theorem \ref{a}. Recall that in general if $G_0 $ is a closed subgroup of $G$ and $N$ is a discrete $G_0 $--module then $H^i (G_0 , N)$ and $H^i (G, \ind ^G _{G_0 }N)$ are naturally isomorphic (see for example \cite{Se1}, I.2.5.) Similarly for the homology of compact $\Lambda (G_0 )$--modules.{\\ \nopagebreak
\hspace*{\fill}$\Box$}
\end{proof}
We would like to thank R.~Sujatha for pointing out how Proposition \ref{o} follows easily from Shapiro's Lemma.

Our reason for including this elementary observation here is because in the case of $G $ isomorphic to $ \z _p ^{\times }$, with $p $ at least $ 3 $, it is closely related to the decomposition of $M$ into eigenspaces for the action of $\Delta \cong { ( \z / p )}^{\times }$, something which is a very useful tool in this situation. In fact, in this case we can specify the 
relationship more precisely.

When $G \cong \z _p^{\times }$ with $p \gge 3$, $G $ is the direct product of 
$G_0 $ and $\Delta $.
Suppose $\chi : \Delta \rightarrow \z _p ^{\times }$ 
is a primitive character in the dual group, $\widehat{\Delta }$, of $\Delta $. 
We denote by $M^{\chi }$ the $\Lambda (G)$--submodule of $M$ picked out by the idempotent of $\z _p [ \Delta ]$ corresponding to $\chi $, and upon which 
$\Delta $ acts via $\chi $. Then there is a direct sum decomposition
\begin{equation}\label{61}
M = \underset{\chi \,\in \, \widehat{\Delta } }{\bigoplus } \, M^{\chi },
\end{equation}
and the $\Lambda (\Gamma )$--rank of $M$ is the sum of the 
$\Lambda (\Gamma )$--ranks of the $M^{\chi }$.

From this we can isolate the 
following stronger result:
\begin{cor}\label{x}
For $G \cong \z _p^{\times }$ and $p \gge 3$
\begin{equation}\label{81}
{\rm{hmrank}}_{\Lambda (G)} \big( M\otimes _{\z _p } \z _p (\chi ) \big) = 
{\rm{rank}}_{\Lambda (\Gamma )} (M^{\chi ^{-1}}),
\end{equation}
where $\Gamma \times \Delta = G$, with $ \Gamma $ and $\Delta $ as above, 
and $\chi $ is any primitive character of 
$\Delta $.
\end{cor}

\begin{Remarks}
\begin{enumerate}
\item We note that it is possible for a finitely 
generated $\Lambda (G)$--module, $M$, to have 
hmrank$_{\Lambda (G)} (M) =0$ but
non-zero $\Lambda (G_0)$--rank, where $G_0$ is any pro-$p$, open subgroup of 
$G$. Indeed this is apparent from the proof of Proposition \ref{o} 
in the special case of $G $ isomorphic to $ \z _p ^{\times }$. Simply take 
$M $ to be $ \Lambda (G)^{\chi }$
for a non-trivial character $\chi : \Delta \rightarrow \z _p^{\times }$.
This is isomorphic to $\Lambda (\Gamma ) (\chi )$, and hence has $\Lambda (\Gamma )$--rank equal to 1.
\item It is also possible to have $\hmrank _{\Lambda (G)} (M)$ negative. 
I thank Otmar Venjakob for pointing out the following example. 
Consider $G \cong \Gamma  \rtimes \Delta $, where $\Gamma $ is taken to be isomorphic to 
$ \z _p $ and 
$\Delta \cong (\z /p )^{\times }$ acts upon $\Gamma $ via a non-trivial 
character, $\omega : \Delta \rightarrow \z _p ^{\times }$. Assume $p $ is at least 
$ 3 $ and 
so $\# \Delta $ is prime to $p$. Then the Hochschild-Serre spectral 
sequence gives $H_i (G,M) = (H_i (\Gamma ,M ))_{\Delta }$ for all $i$. 
Take $M \cong \z _p (\omega )$. Then $H_0 (G,M) $ vanishes and 
$H_1 (G,M) = (\Hm (\Gamma ,\z _p (\omega )))_{\Delta } =\z _p $.
Thus $\hmrank _{\Lambda (G)}(\z _p (\omega )) =-1$.
In this case $G$ is a Poincar\'{e} group with dualising module 
$\q _p / \z _p ( \omega )$, and I have been unable to find such an 
example where $G$ has trivial dualising module.
\end{enumerate}
\end{Remarks}

\begin{question}
If $G$ is a compact, open subgroup with no $p$--torsion in a connected, 
reductive algebraic group over $\q _p $ then $G$ is a Poincar\'{e} group 
whose dualising module is $\q _p /\z _p $ with trivial $G$--action. (This is 
essentially due to Lazard, \cite{La65} with Burt Totaro pointing out how 
to extend the argument to $G$ not necessarily a pro-$p$ group.) In this 
case, is it true that the $ \hmrank _{\Lambda (G)}(M) $ is always positive 
for all finitely generated $\Lambda (G)$--modules?
\end{question}

We would hope that these ideas would generalise, at least to the case 
$G \cong G_0 \rtimes \Delta$, but there are added complications, 
especially if $\Delta $ is non-Abelian or the representations of 
$\Delta $ do not decompose completely into characters over $\z _p $. 
If we cannot exhibit $\Delta $ as a subgroup of $G$ then it is 
presently not at all clear how to proceed.

\subsection{{\it{\bfseries{M}}} fails to be finitely generated}\label{2.2}

A second situation in which the homological definition of 
rank is of some interest 
is in considering modules which are not finitely generated.
The general situation is rather unclear to us and so, instead of
giving a general discussion, we will consider a particular example of
one situation which arose in \cite{CoHo99}. Let $E$ be an 
elliptic curve defined 
over a number field, $F$, and with no complex multiplication over the 
algebraic closure of $F$. Denote by 
$\E $ the set of all $p^{\rm{th}}$--power torsion points in $E$, and by 
$\F $ we mean $ F(\E )$, the field of definition of $\E $. Let $\G $ denote 
the Galois group of $\F $ over $F$.
It is well known, due to Serre, that the action of
$\G $ on $\E $ gives an injection of $\G $ into GL$_2(\z _p )$
as an {\em{open}} subgroup, under our assumption that $E$ has no complex 
multiplication. For $p $ at least $5$, the group GL$_2 (\z _p )$ contains
no element of order $p$ and thus, from the properties of $p$-adic Lie groups
listed above, $\G $ has $p$--cohomological dimension equal to 4. We denote 
by $F_n $ the field $F(E_{p^{n+1}} )$ obtained by adjoining the $p^{n+1}$
torsion points on $E$ to $F$, and by $F^{cyc}$ the cyclotomic 
$\z _p $--extension of $F$. As is well known, it follows from the Weil 
pairing that $F^{cyc}$ is contained in $\F $. We define
\begin{equation}\label{32}
H = \Gal (\F /F^{cyc} ),\;\;H_0 = \Gal (\F /F_0^{cyc} ),\;\;
{\rm{and}}\;\;\Gamma = 
\Gal (F^{cyc} /F ).
\end{equation}
Let $S$ be any finite set of primes of $F$ containing all primes dividing 
$p$, Archimedean primes of $F$ and primes at which $E$ has bad reduction,
and let $\fs $ denote, as before, the maximal extension of $F$ unramified
at all primes not contained in $S$.
Then, since $\F $ is an unramified extension of $F$ at primes not in $S$, 
$\F $ is contained in $\fs $.
For any, not necessarily finite, extension $L$ of $F$ contained in 
$\F $ the $p^{\infty }$--Selmer 
group of $E$ over $L$, which we denote by $\Sel _p (E/L)$, can be defined by 
the exactness of the sequence
\begin{equation}\label{33}
0 \longrightarrow \Sel _p (E/L) \longrightarrow 
H^1(F^S /L ,E_{p^{\infty }} ) \overset{\lambda _L}{\longrightarrow}
\underset{\nu \in S}{\bigoplus}\, J_{\nu } (E/L).
\end{equation}
(More usually local conditions are imposed at all primes, but the choice of the primes in $S$ ensures that this definition is equivalent.)
Here we are using the notation 
\begin{equation}\label{34}
J_{\nu } (E/L) = \underset{\longrightarrow}{\lim}\,\underset{\omega \mid \nu 
}{\bigoplus}H^1(K_{\omega } ,E)(p),
\end{equation}
where $K$ runs over the finite extensions of $F$ contained in $L$, 
$\omega $ denotes a prime of $K$ dividing $\nu $ and the direct limit is 
taken with respect to the restriction maps.

Then $\G $ acts continuously on $\Sel _p (E/\F )$, where the latter module 
is given the discrete topology, thus this action extends to a continuous 
action of the Iwasawa algebra, $\Lambda (\G )$. We will in fact mainly 
consider the Pontrjagin dual 
\begin{equation}\label{35}
\CC _p (E/\F ) = \widehat{\Sel _p (E/\F )}
\end{equation}
which has the structure of a compact $\Lambda (\G )$--module. Similarly, 
$\CC _p (E/F^{cyc})$ is a compact $\Lambda (\Gamma )$--module. It is well 
known that $\CC _p (E/F^{cyc})$ is a finitely generated 
$\Lambda (\Gamma )$--module. Similarly $\CC _p (E/\F )$ is a 
finitely generated 
$\Lambda (\G )$--module (see \cite{CoHo99}, Theorem 2.9.)

The following conjecture is folklore.
\begin{conj}\label{i}
Assume $p \gge 3$. If $L$ is any extension of $F^{cyc}$ contained in $F^S$ 
then 
\begin{enumerate}
\item the map $\lambda _L$ in the exact sequence (\ref{33}) 
is a surjection and 
\item the group $H^2(F^S/L,E_{p^{\infty }})=0$.
\end{enumerate}
\end{conj}
In fact, as has been pointed out to us by Y.~Ochi, it is known that 
\begin{equation}\label{43}
H^2(F^S/\F ,\E )=0,
\end{equation}
but this is still only a conjecture in general for arbitrary extensions $L$ 
of $F^{cyc}$. We have restricted to $p \gge 3$ because the vanishing of 
$H^2(F^S/L,E_{p^{\infty}})$ can fail for $p=2$.

In the case $L=F^{cyc}$ this is equivalent to a longstanding conjecture on 
the rank of $\CC _p (E/F^{cyc})$ as a $\Lambda (\Gamma )$--module, 
originally due to Mazur in the case $E$ has good, ordinary reduction at 
all primes of $F$ dividing $p$. In particular, if $E$ has potential good 
ordinary reduction or potential multiplicative reduction at all primes of 
$F$ dividing $p$ then for $F^{cyc} $ 
the following is equivalent to conjecture \ref{i}.
\begin{conj}\label{j}
$\CC _p (E/F^{cyc} )$ is $\Lambda (\Gamma )$--torsion.
\end{conj}
It follows from the structure theorem for finitely generated 
$\Lambda (\Gamma )$--modules that $\CC _p (E/F^{cyc} )$ has finite 
$\z _p $--rank in this case (though it will generally not be finitely 
generated as a $\z _p $--module.) We obtain Conjecture \ref{i} for
general $L$ from the 
conjecture for $F^{cyc}$ by taking the direct limit over {\em{all }} finite 
extensions of $F$ contained in $F^S$.

It was noted in \cite{CoHo99} that the classical property of being torsion 
and having $\mu $-invariant zero in the theory of $\z _p $--extensions seems 
to have a natural parallel in the GL$_2$ theory of the extension $\F $ of $F$.
We recall that for $\Gamma$ any group isomorphic to $\z _p $, if $M$ is a 
finitely generated $\Lambda (\Gamma )$--module then $M$ is torsion and 
has  
$\mu $-invariant zero if $M$ is actually finitely generated over $\z _p $.
This is an immediate consequence of the definition of the classical Iwasawa 
$\mu $-invariant.

We have given above in (\ref{203}) a definition of $\mu $-invariant for 
general $\Lambda (G)$--modules when $G$ is a pro-$p$, $p$-adic, Lie group. 
But certain results discussed in \S 6 of \cite{CoHo99} for Selmer groups of 
elliptic curves and in \cite{Ho1} for other situations indicate that another 
analogue of the $\mu $-invariant being zero in this GL$_2 $ situation is that 
a finitely generated $\Lambda (\G )$--module, $M$, be in fact finitely 
generated over the subalgebra, $\Lambda (H)$. Unfortunately this is 
stronger, the two notions are not equivalent.
\begin{lem}\label{yy}
Let $M$ be a $\Lambda (\G )$--module which is finitely generated over $\Lambda (H)$. Then $\mu (M) =0 $ where $\mu $ here is as defined in (\ref{203}), for any choice of pro-$p$, open subgroup $G' $ of $\G $.
\end{lem}
\begin{proof}
The $p$-torsion submodule $M(p)$ is a finitely generated $\Lambda (H)$--submodule, and so all of the sub-quotients, $p^i M(p) / p^{i+1} $, have rank 0 over $\Omega (G')$.{\\ \nopagebreak
\hspace*{\fill}$\Box$}
\end{proof}
The converse clearly fails. Consider, for example, the module 
$ \Lambda (\Gamma ) /p $ with the structure of a $\Lambda (\G )$--module via the natural projection $\Lambda (\G ) \twoheadrightarrow \Lambda (\Gamma )$. This example, however, is rather small. O.~Venjakob has given a definition of pseudonull modules for $\Lambda (G)$--modules which extends the usual 
definition in the case of $G \cong \z _p $, \cite{Ven}. This example is a pseudonull $\Lambda (\G )$--module.

In a similar vein, the $\Lambda (H_0)$--rank of 
$M$, as defined in (\ref{4}) above, bears similarities 
to the $\lambda $-invariant in the theory 
of $\z _p $--extensions. (One has to restrict to $H_0$ for the usual reasons, 
to ensure that $\Lambda (H_0 )$ has no non-trivial torsion.) We 
would like to propose 
that a more suitable invariant should be the homological $\Lambda (H)$--rank 
defined above. As we shall see, this has computational advantages for modules which are not finitely generated over $\Lambda (H_0 )$.

\begin{notn}
In the sequel we let $r$ denote the number of primes of $F^{cyc}$ at which $E$ has split multiplicative reduction.
\end{notn}
\begin{thm}\label{k}
Assume $p \gge 5$ and $E$ does not have potential good supersingular 
reduction at any prime of $F$ dividing $p$. Assume conjecture \ref{i} 
for the fields
$\F $ and $F^{cyc}$. Then $\CC _p (E/\F )$ has finite homological rank 
over $\Lambda (H)$. If $F$ contains $\mu _p $, the $p^{th}$ roots 
of unity, then this is given by
\begin{equation}\label{36}
{\text{\rm{hmrank}}}_{\Lambda (H)} \big( \CC _p (E/\F ) \big) = \rk _{\z _p }
\big( \CC _p (E/F^{cyc} ) \big) + r.
\end{equation}
\end{thm}
\begin{proof} We have the following fundamental diagram, with vertical maps given by restriction maps and exact rows as a consequence of our assumption concerning Conjecture \ref{i}:
\begin{equation}\label{40}
\begin{array}{cccccccc}
0 \rightarrow &\Sel _p (E/\F )^H &\rightarrow &H^1(F^S/\F ,\E )^H 
&\overset{\phi _{\infty }}{\longrightarrow} & \big( \underset{\nu \in S}{\bigoplus}\, J_{\nu } (E/\F )  \big) ^H  & {} & {} \\
{} & \uparrow {\scriptstyle{f}} & {} & \uparrow {\scriptstyle{g}} & { } &
\uparrow {\scriptstyle{h}} & {} & {} \\
0 \rightarrow & \Sel _p (E/F^{cyc}) & \rightarrow & H^1(F^S/F^{cyc} ,\E ) 
&\overset{\lambda _{F^{cyc}}}{\longrightarrow} &  \underset{\nu \in S}{\bigoplus}\, J_{\nu } (E/ F^{cyc} ) & \rightarrow & 0
\end{array}
\end{equation}
This diagram is analysed in detail in \cite{CoHo99}, where
the following is proved (Lemma 6.7.)
\begin{lem}\label{l}
\begin{enumerate}
\item The kernel and cokernel of $g$ are both finite.
\item The cokernel of $h$ is finite and the kernel of $h$ has finite 
$\z _p $--corank, at most equal to $r$.
\item This is an equality if $F$ contains $\mu _p $.
\end{enumerate}
\end{lem}
Since $\CC _p (E/F^{cyc})$ has finite 
$\z _p $--rank, by assumption, it follows from the snake lemma that 
$\CC _p (E/\F )_H $ has finite $\z _p $--rank, given by
\begin{equation}\label{41}
\begin{array}{ll}
\rk _{\z _p } \big( \CC _p (E/\F )_H \big) & = \rk _{\z _p } 
\big( \CC _p (E/F^{cyc}) \big) + \rk _{\z _p } \big( \widehat{\kr (h)} \big) \\
{} & \lle \rk _{\z _p } \big( \CC _p (E/F^{cyc} ) \big) + r,
\end{array}
\end{equation}
with equality if $F$ contains $\mu _p $.

Since we are assuming Conjecture \ref{i} for $\F $, the top row of 
(\ref{40}) extends to a long exact sequence in cohomology:
\begin{equation}\label{42}
\begin{array}{c}
0 \rightarrow {\rm{Coker}} (\phi _{\infty }) \rightarrow H^1 \big( H ,\Sel _p 
(E/\F ) 
\big) \rightarrow \cdots \\
\rightarrow 
H^3 \big( H , H^1(F^S/\F , \E ) \big) \rightarrow
H^3 \big( H , \underset{\nu \in S}{\bigoplus}\, J_{\nu } (E/\F )  \big) 
\rightarrow 0
\end{array}
\end{equation}
Since the map $\lambda _{F^{cyc}}$ is surjective and ${\rm{Coker}} (h)$ is finite, it
follows from diagram (\ref{40}) and the snake lemma that the cokernel of 
$\phi _{\infty }$ is finite.
\begin{lem} \label{m}
We assume $p \not= 2$. Then it follows from our assumption of the vanishing of
$H^2(F^S/F^{cyc},\E )$ that the groups 
$H^i \big( H,H^1(F^S/\F,\E ) \big) $ are finite for all $i \gge 1$.
\end{lem}
\begin{proof}
We exploit the fact that it is known that $H^j(F^S/\F ,\E ) =0$ for 
$j $ at least $ 2$. 
Then it follows from the Hochschild-Serre spectral sequence 
that we have the exact sequence
\begin{equation}\label{45}
H^i(F^S/F^{cyc},\E ) \longrightarrow H^{i-1} \big( H, H^1(F^S/\F ,\E ) \big) 
\longrightarrow H^{i+1}(H,\E ),
\end{equation}
for $i $ at least $ 2$. 
We have assumed that $H^2(F^S/F^{cyc},\E )$ vanishes, whilst 
for $i \gge 3 $ the groups $H^i(F^S/F^{cyc},\E )$ vanish since 
the $p$ cohomological dimension of $\Gal (F^S/F^{cyc})$ equals $2$ under our assumption that $p $ is not equal to $2$.
The groups $H^i(H,\E )$ 
are easily seen to be finite for all $i \gge 0$. See, for example, 
\cite{CoSu99}.{\\ \nopagebreak
\hspace*{\fill}$\Box$}
\end{proof}
\begin{lem} \label{n}
Assume that $p \gge 5$. Then for all $i \gge 1$ one has
\begin{equation}\label{46}
H^i(H,\underset{\nu \in S}{\bigoplus}\, J_{\nu } (E/\F ))=0.
\end{equation}
\end{lem}
\begin{proof}
It is proven in \cite{CoHo99} that (\ref{46}) is true with $H$ replaced by 
$\G $ (Proposition 5.12 for $\nu \nmid p$ and Corollary 5.23 for $\nu \mid p$.)
Since $cd _p (\Gamma )=1$ the Hochschild-Serre spectral sequence gives rise 
to short exact sequences
\begin{equation}\label{47}
\begin{array}{c}
0 \rightarrow H^j \big( \Gamma , H^j(H,\underset{\nu \in S}{\bigoplus}\, 
J_{\nu } (E/\F ) )   \big) \rightarrow   \\
H^{j+1} \big(   \G , 
\underset{\nu \in S}{\bigoplus}\, J_{\nu } (E/\F ) \big)   
\rightarrow H^{j+1}  \big( H, \underset{\nu \in S}{\bigoplus}\, 
J_{\nu } (E/\F )  \big) ^{\Gamma }  \rightarrow 0
\end{array}
\end{equation}
Since $H^{j+1} (H,\underset{\nu \in S}{\bigoplus}\, J_{\nu } (E/\F )  ) $ is
$p$--primary and 
has the discrete topology, and since $\Gamma $ is isomorphic to $ \z _p $, 
the vanishing of 
$H^{j+1} (H,\underset{\nu \in S}{\bigoplus}\, J_{\nu } (E/\F )  ) ^{\Gamma} $
for $j \gge 0$
implied by (\ref{47}) ensures that $H^i(H,\underset{\nu \in S}{\bigoplus}\, J_{\nu } (E/\F ) )$ vanishes for all $i \gge 1$.{\\ \nopagebreak
\hspace*{\fill}$\Box$}
\end{proof}
Thus we see from the exact sequence (\ref{42}), Lemmas \ref{m} and \ref{n}
and the finiteness of $\cok (\phi _{\infty })$ that for all $i \gge 1$ the 
$H^i(H,\Sel _p (E/\F ))$ are finite.

The duality between homology and cohomology described in 
(\ref{24}) extends to 
this situation, even though $\CC _p (E/\F )$ is not finitely generated as 
a $\Lambda (H)$--module. 
Indeed, taking Pontrjagin duals passes between the categories of compact $\Lambda (H)$--modules and discrete $\Lambda (H)$--modules, without any finite generation assumption. This duality again sends compact, projective modules to discrete, injective modules, and we still have 
\begin{equation}\label{48}
\widehat{M_H} = \widehat{M} ^H ,
\end{equation}
for any compact $\Lambda (H)$--module, without assuming $M$ to be finitely generated. In fact, since finitely generated $\Lambda (G)$--modules are pseudo-compact as $\Lambda (H)$--modules, this is explained in greater detail in \cite{Br66}. Thus we also still have a canonical isomorphism
\begin{equation}\label{49}
H^i(H,\widehat{M}) \cong \widehat{H_i (H,M)}.
\end{equation}

It follows that $H_i (H,\CC _p (E/\F ))$ 
is finite for all $i \gge 1$, and so the 
homological $\Lambda (H)$--rank of $\CC _p (E/\F )$ is just the $\z _p$--rank
of $\CC _p (E/\F )_H$, given by (\ref{41}). This completes the proof of Theorem \ref{k}{\\ \nopagebreak
\hspace*{\fill}$\Box$ }
\end{proof}

As a consequence of Theorem \ref{a}, in the case $H$ is pro-$p$ and 
$\CC _p (E/\F )$ is actually finitely generated over $\Lambda (H)$
we have calculated the true rank, not just the homological rank.
Thus we obtain directly the following result (Proposition 6.9) in 
\cite{CoHo99}.

\begin{cor}\label{xx} Assume the hypotheses of Theorem \ref{k}. Assume also 
that $H$ is a pro-$p$ group and that $\CC _p (E/F^{cyc} )$ 
has $\mu $-invariant equal to 0.
Then if $L$ is any finite extension of $F^{cyc}$ contained in $\F $ 
\begin{equation}
\lambda _{inv} \big( \CC _p (E/L ) \big) = \mid L:F^{cyc} \mid
\Big( \lambda _{inv} \big( \CC _p (E/F^{cyc} ) + r   \big)    \Big) - r_L
\end{equation}
where $r_L $ denotes the number of primes of $L$ at which $E$ has split 
multiplicative reduction.
\end{cor}

This is clear from the fact that if $M$ is a finitely generated 
$\Lambda (H)$--module of rank $d$ and if $H'$ is an open subgroup of 
$H$, then $M$ is a finitely generated $\Lambda (H')$--module of 
rank $\mid H : H' \mid$. Corollary \ref{xx}  
was originally obtained as a consequence of Hachimori and Matsuno's 
more general formula concerning the growth of $\lambda $-invariants in 
a $p$--extension (see \cite{HaMa98}). We believe the 
above proof explains, in this situation, 
the origin of the terms in their formula.

We conclude this section by illustrating Theorem \ref{k} with an example expanding upon 
one considered in \cite{CoHo99}.

\begin{ex} Take $E$ to be the elliptic curve $X_0(11)$ with minimal 
Weierstra{\ss} equation over $\q $ given by
\begin{equation}\label{70}
E:y^2 + y = x^3 -x^2 -10x -20.
\end{equation}Take $p=5$, a prime at which $E$ has good, ordinary 
reduction. Let $H$ be the group 
$\Gal (\q (E_{5^{\infty }}) /\q (\mu _{5^{\infty }}) )$.
It follows from the work of Lang and Trotter \cite{LaTr} that $H$ is pro-$p$.
Then it is explained in \S 6 and 7 of \cite{CoHo99} that
$\CC _5 (E/\q (E_{5^{\infty }}) )$ is not finitely generated as a 
$\Lambda (H)$--module. However, Corollary 7.10 of \cite{CoHo99} states
\begin{equation}\label{71}
{\rm{dim}}_{K(H)} \big( \,  \CC _5 (E/\q 
(E_{5^{\infty }}))\otimes_{\Lambda (H)} K(H) \,  \big) =4.
\end{equation}
Recall that $K(H)$ denotes the skew field of fractions of $\Lambda (H)$. 
It follows that $\CC _5 (E/\q (E_{5^{\infty }}) )$ actually has finite 
$\Lambda (H)$--rank.

The only rational prime at which $E$ has non-integral $j$-invariant is 11. 
There are 4 primes of $\q (\mu _{5^{\infty }} )$ over 11. Thus the number $r$ 
appearing in Theorem \ref{k} is 4. 
It follows from the isogeny invariance of the $\lambda $-invariant and the 
results in \cite{CoHo99} that $\CC _5 (E/\q (\mu _5^{\infty } ) )$ has 
$\z _5 $--rank equal to 0 (although it has non-zero $\mu $-invariant and 
so is not finite.) Thus Theorem \ref{k} in this case states:
\begin{equation}\label{72}
{\rm{hmrank}}_{\Lambda (H) } \big( \, \CC _5 (E/\q (E_{5^{\infty }}) ) 
\, \big) =4,
\end{equation}
agreeing with the actual $\Lambda (H)$--rank in (\ref{71}).
\end{ex}

We have unfortunately so far 
been unable to deduce from Theorem \ref{k} the answer 
to the original question raised in \cite{CoHo99}, of whether the 
$\Lambda (H_0 )$--rank is always finite under these conditions. As illustrated
by the above, in the very few examples discussed in 
\cite{CoHo99} where we 
know the $\Lambda (H_0 )$--rank of $\CC _p (E/\F )$ it agrees with the 
$\Lambda (H_0 )$--hmrank.

\section{Variation of the ${\mathbb{\mu }}$-invariant of Elliptic Curves under Isogeny}

Recalling the definition (\ref{203}) of the $\mu $-invariant of modules finitely generated over Iwasawa algebras given in \S \ref{1.2}, we consider how this behaves for Selmer groups of elliptic curves. Perrin-Riou in \cite{P-R89} and the appendix of \cite{P-R87} and Schneider in \cite{Schn87} consider the variation under isogeny of the classical Iwasawa $\mu $-invariant of finite dimensional $p$-adic representations of $\Gal ( {\overline{\q }} / \q)$. For simplicity, we will restrict to Selmer groups of elliptic curves, but as in \cite{P-R89}, the argument holds more generally if we use Greenberg's definition of Selmer groups for $p$-adic representations given in \cite{Gr89} and make the appropriate generalisation of the assumption {\bf{L}} below. We continue with the situation under discussion in \S \ref{2.2} and maintain all the same notation. 

The behaviour of the $\mu $-invariant for Selmer groups of elliptic curves under isogeny has also been considered independently by J.H.~Coates and R.~Sujatha, who also have an alternative approach to these results.

Assume $E_1 $ and $E_2 $ are non-CM elliptic curves defined over a number 
field, $F$, and they are isogenous via  $\phi : E_1 \rightarrow E_2 $. We 
suppose $\F = F( E_{1,p^{\infty}})$ is a pro-$p$ extension of $F$. Note 
this field is the same as $F(E_{2,p^{\infty }}) $, see
\cite{CoHo99} Lemma 7.8. Let $A $ denote the $p$-part of the kernel of 
$ \phi $. 
Throughout this section we need to assume a special case of conjecture \ref{i}.
\begin{ass}{\bf{L}} The localisation maps
\begin{equation}
\lambda _{i,\F } : H^1(\fs /\F ,E_{i,p^{\infty }}) \longrightarrow \underset{\nu \in S}{\bigoplus} J_{\nu }(E_i /\F ) 
\end{equation}
are surjective for both $i =1, 2$.
\end{ass}
In fact, as explained in \cite{CoHo99} Proposition 7.9, if this holds 
for one of the curves $E_i $ then it automatically holds for the other.

The isogeny, $\phi $, induces homomorphisms on the cohomology groups and Selmer groups:

\begin{equation}\label{213}
\begin{array}{ccccccc}
0 \longrightarrow & \Sel _p (E_2 /\F ) & \longrightarrow & H^1 (\fs /\F , E_{2 , p^{\infty}} ) & \overset{\lambda _{2,\F }}{\longrightarrow} & \underset{\nu \mid S}{\bigoplus} J _{\nu } (E_2 /\F ) & \longrightarrow 0 \\
{} & \uparrow \phi _1 & {} & \uparrow \phi _2 & {} &
\uparrow \phi _3 & {} \\
0 \longrightarrow & \Sel _p (E_1 /\F ) & \longrightarrow & H^1 (\fs /\F , E_{1 , p^{\infty}} ) & \overset{\lambda _{1,\F }}{\longrightarrow} & \underset{\nu \mid S}{\bigoplus} J _{\nu } (E_1 /\F ) & \longrightarrow 0 \\
\end{array}
\end{equation}

We wish to consider the change in the $\mu $-invariant.

\begin{thm}\label{aa}
Let $p \gge 5$. Suppose $\phi : E_1 \longrightarrow E_2$ is an isogeny between two elliptic curves which do not admit complex multiplication, as described above. Assume also that the $E_i $ satisfy Assumption {\bf{L}} and that $\F $ is a pro-$p$ extension of $F$. Then, as $\Lambda (\G )$--modules
\begin{equation}\label{212}
\begin{array}{c}
\mu (\widehat{\cok (\phi _1 )}) - \mu (\widehat{\kr (\phi _1 )}) = \underset{\nu \mid \infty }{\sum } \,  \ord _p \big( \# (A(F_{\nu }) ) \big) \;  \\
- \; \mid F:\q \mid \ord _p (\# A) \; - \;
\underset{\nu \mid p }{\sum } \, \ord _p \big( \mid \# \tilde{A} _{\nu } \mid _{\nu } \big)
\end{array}
\end{equation}
where $A = \kr (\phi (p))$ as above, $\nu $ denotes a prime of $F$, by $A(F_{\nu })$ we mean the set of $F_{\nu }$-valued points in $A$, and 
$\tilde{A} _{\nu } $ denotes the image of $A$ under the reduction map at $ \nu \,: \, E_1 \longrightarrow \tilde{E_1} $.
\end{thm}
\begin{proof}
This follows closely the ideas of Perrin-Riou in the cyclotomic case given in detail for elliptic curves in the appendix to \cite{P-R87}. We proceed similarly to the proof of Theorem \ref{k} above, with a close analysis of a diagram, this time (\ref{213}). 

The snake lemma gives 
\begin{equation}\label{214}
\begin{array}{c}
0 \longrightarrow \kr (\phi _1 ) \longrightarrow \kr (\phi _2 ) \longrightarrow \kr (\phi _3 ) \longrightarrow \cok (\phi _1 ) \\ \longrightarrow 
\cok (\phi _2 ) \longrightarrow \cok (\phi _3 ) \longrightarrow 0.
\end{array}
\end{equation}
The first thing to remark is that the Pontrjagin dual of the kernels and cokernels of $\phi _2 $ and $\phi _3 $ (and hence of $\phi _1 $ by the snake lemma) are all finitely generated $\Lambda (\G )$--modules, annihilated by a finite power of $p$. 
Note that although Corollary \ref{z} is given in terms of group homology, because of (\ref{24}) it can also be rephrased in terms of cohomology as for Theorem \ref{a}. For discrete, cofinitely generated $\Lambda (\G) $--modules, we will use also the notation
$\chi (\G ,D) $ to denote $ \underset{i \gge 0}{\prod} \# \big( H^i(\G ,D) \big) ^{(-1)^i } $ where the choice between cohomology and homology is dictated by whether the module is discrete or compact, respectively. Since either choice gives zero for finite modules when $\G $ is pro-$p$, there will be no ambiguity. Thus if $M$ is a finitely generated $\Lambda (\G )$--module then $\mu (M) = \ord _p (\chi (\G , \widehat{M(p)} ))$ also.
Then the multiplicativity of $\chi (\G ,\_ ) $ along exact sequences gives
\begin{equation}\label{215}
 \frac{\chi (\G , \cok (\phi _1)) }{ \chi (\G , \kr (\phi _1 ))} =
\frac{\chi (\G , \cok (\phi _2))}{ \chi (\G , \kr (\phi _2)) } \times 
\frac{\chi (\G , \kr (\phi _3 )) }{ \chi (\G , \cok (\phi _3 ))},
\end{equation}
where all Euler characteristics are defined. Since $\cok (\phi _1)$ and $\kr (\phi _1 )$ are both annihilated by a finite power of $p$, by Corollary \ref{dd} the number we require for the Theorem is then the logarithm to base $p$ of this.

For the first term, recall that we know $H^i (\fs /\F , E_{i, p^{\infty }})$ vanishes for $i \gge 2$, by (\ref{43}) and the assumption that $p \gge 3$. Thus by the long exact sequence in cohomology of which the central vertical arrow is a part, we have $\cok (\phi _2 ) = H^2 (\fs /\F ,A)$. Since isogenies are surjective $H^0 (\fs /\F ,E_{1,p^{\infty }} ) = E_{1,p^{\infty }}$, surjects onto $H^0(\fs /\F , E_{2,p^{\infty }}) = E_{2,p^{\infty }}$, and so 
$\kr (\phi _2 ) = H^1 (\fs /\F ,A)$. The full Hochschild-Serre spectral sequence:
\begin{equation}
H^i (\G , H^j (\fs /\F , A )) \Rightarrow H^{i+j} (\fs /F ,A)  
\end{equation}
is bounded since $\G $, $\Gal (\fs /\F )$ and $\Gal (\fs /F )$ all have finite cohomological dimension at $p$. It follows from this that
\begin{equation}\label{216}
\frac{\chi (\G , H^0 (\fs /\F ,A)) \chi (\G , H^2 (\fs /\F ,A)) }{ \chi (\G , H^1 (\fs /\F ,A)) } =  \chi (\Gal (\fs /F ),A ).
\end{equation}
Precise details of how Euler Characteristics behave in spectral sequences, and in particular how to obtain this formula, are given in \cite{Ho1}. 
But since $H^0 (\fs /\F ,A) $ is just $A$ and in particular is finite, as remarked in \S \ref{1.2} it has Euler characteristic equal to one, and thus 
\begin{equation}\label{217}
\frac{\chi (\G ,\cok (\phi _2 ))}{\chi (\G , \kr (\phi _2 ))} = \chi (\Gal (\fs /F ),A ).
\end{equation}
Locally, for any prime $\nu $ of $F$ in $S$ we define
\begin{equation}
\begin{array}{lll}
C_{\nu }^i (\F ) = & \underset{\longrightarrow}{\lim} \, \underset{\omega \mid \nu }{\prod} H^i (F_{n,\omega},A) & {\rm{if}}\;\, \nu \nmid p \\
{} & \underset{\longrightarrow}{\lim} \, \underset{\omega \mid \nu }{\prod} H^i (F_{n,\omega}, \tilde{A} _{\nu }) & {\rm{if}}\;\, \nu \mid p
\end{array}
\end{equation}
where the direct limits are taken with respect to the restriction homomorphisms. Recall that $\tilde{A} _{\nu } $ is the image of $A$ under the reduction map at $ \nu \, :\, \,E_1 \longrightarrow \tilde{E_1} $. As explained in the local calculations in \cite{CoHo99} \S 5, it follows from the results in \cite{CoGr96} that 
\begin{equation}
\begin{array}{lll}
J_{\nu } (E_i /\F ) = & \underset{\longrightarrow}{\lim}\,\underset{\omega \mid \nu }{\prod } H^1(F_{n,\omega },E_{i,p^{\infty }}) & {\rm{if}}\;\, \nu \nmid p \\
{} & \underset{\longrightarrow}{\lim}\,\underset{\omega \mid \nu }{\prod } H^1(F_{n,\omega },\tilde{E} _{i,p^{\infty }}) & {\rm{if}}\;\, \nu \mid p. 
\end{array}
\end{equation}
Let $X $ be $ E_{i,p^{\infty }}$ if $\nu \nmid p$, whilst if $\nu \mid p$ then $X$ denotes $\tilde{E} _{i,p^{\infty }} $.
Then the groups  $\underset{\longrightarrow}{\lim}\,\underset{\omega \mid \nu }{\prod } H^2(F_{n,\omega },X)$ vanish (also explained in \cite{CoHo99} \S 5) and as described above for the global map $\phi _2 $, we have the simple descriptions  $\kr (\phi _3 ) = \underset{\nu \in S}{\bigoplus } \, C^1_{\nu } (\F )$ and $\cok (\phi _3 ) = \underset{\nu \in S}{\bigoplus } \, C^2_{\nu } (\F )$.

For each prime $\nu $ in $S$ we fix a choice of prime $\F $ lying above, which will also be denoted by $\nu $. It follows from Shapiro's Lemma that
\begin{equation}
H^j (\G , C^i_{\nu } (\F )) = H^j (\Delta _{\nu } , H^i (F_{\infty , \nu},X)),
\end{equation}
where $X$ is as above and $\Delta _{\nu }$ is the decomposition group of $\nu $ in $\G $, that is the Galois group of the extension $F_{\infty , \nu} / F_{\nu }$. Thus
\begin{equation}
\chi (\G ,\kr (\phi _3 )) = \underset{ \nu \in S , \, \nu \mid p}{\prod} \,
\chi (\Delta _{\nu } ,H^1( F_{\infty , \nu} , \tilde{A} _{\nu } )) \times 
\underset{ \nu \in S , \, \nu \nmid p}{\prod} \,
\chi (\Delta _{\nu } ,H^1( F_{\infty , \nu} , {A} ))
\end{equation}
and similarly for $\cok (\phi _3 )$. Then, as we saw when considering $\phi _2 $, the Hochschild-Serre spectral sequence
\begin{equation}
H^i (\Delta _{\nu } ,H^j ( F_{\infty , \nu} , X)) \Rightarrow H^{i+j} (F_{\nu } ,X)
\end{equation}
gives 
\begin{equation}\label{219}
\frac{\chi (\G ,\kr (\phi _3 ))}{\chi (\G ,\cok (\phi _3 ))} = \frac{1}{\underset{ \nu \in S , \, \nu \mid p}{\prod} \, \chi (\Gal ( \overline{F_{\nu }}/F_{\nu }),\tilde{A} _{\nu }) \times 
\underset{ \nu \in S , \, \nu \nmid p}{\prod} \, \chi (\Gal ( \overline{F_{\nu }}/F_{\nu }),{A}) }.
\end{equation}
Since $A$ is a $p$--group, it is known, due to Tate \cite{Mi86} Lemma I.2.9, that the Euler Characteristic corresponding to primes in $S$ which do not lie above $p$ equals one. Then substituting (\ref{217}) and (\ref{219}) into (\ref{215}) and taking the logarithm to base $p$ gives
\begin{equation}\label{402}
\begin{array}{c}
\mu \big( \widehat{\cok ( \phi _1 )}   \big) - \mu \big( \widehat{\kr ( \phi _1 )}   \big) =
\ord _p \big( \chi ( \Gal (\fs /F ),A  )   \big)  \\
- \underset{\nu \in S, \, \nu \mid p}{\sum } \, 
\ord _p \big( \chi ( \Gal ( \overline{F_{\nu }} / F_{\nu }) , \tilde{A} _{\nu }    )   \big).
\end{array}
\end{equation}
We recall here Tate's formulae for calculating such Euler characteristics. If $M$ is a finite $\Gal (\fs /F )$--module such that all primes $\nu $ of $F$ dividing the order of $M$ are contained in $S$, then (see \cite{Mi86} Theorem I.5.1)
\begin{equation}\label{400}
\chi (\Gal (\fs /F ),M) = \underset{\nu {\text{ Arch}}}{\prod} \, \frac{\# H^0 ({\overline{F _{\nu }}} /F_{\nu },M)}{\mid \# M \mid _{\nu }}.
\end{equation}
If $\nu $ is a non-Archimedean prime of $F$ and $M$ is a finite $\Gal ({\overline{F_{\nu }}}/F_{\nu })$--module, then (see \cite{Mi86} Theorem I.2.8)
\begin{equation}\label{401}
\chi (\Gal ({\overline{F_{\nu }}}/F_{\nu }),M) = \mid \# M \mid _{\nu }.
\end{equation}
Substituting these into (\ref{402}) gives precisely the formula (\ref{212}) in the Theorem.{\\ \nopagebreak
\hspace*{\fill}$\Box$}
\end{proof}
\begin{cor}\label{tt}
Maintaining the assumptions of Theorem \ref{aa}, suppose also that the Pontrjagin duals of the $p^{\infty }$--Selmer groups, $\CC _p(E_i /\F )$, are both $\Lambda (\G )$--torsion. Let $L$ be any finite extension of $F$ contained in $\F $. Let $G_L $ denote $ \Gal (\F /L )$, an open subgroup of $\G $. If $M$ is a finitely generated $\Lambda (\G )$--module, we use the notation $\mu _L $ to denote the $\mu $-invariant of $M$, as defined by (\ref{203}), as a $\Lambda (G_L )$--module. Then
\begin{equation}\label{405}
\begin{array}{c}
\mu _L (\CC _p (E_2/\F )) - \mu _L (\CC _p (E_1/\F )) =
\underset{\omega \mid \infty}{\sum} \, \ord _p \big( \# A (L_{\omega }) \big) \; \\
- \; \mid L : \q \mid \ord _p ( \# A ) \; - \; 
\underset{\omega \mid p}{\sum} \, \ord _p \big( \mid \# \tilde{A} _{\omega } \mid _{\omega } \big) 
\end{array}
\end{equation}
where $\omega $ denotes a prime of $L$, and other terms are then as defined in Theorem \ref{aa}.
\end{cor}
\begin{Remark}
We note that, as explained in \cite{CoHo99} Proposition 7.9, the $\Lambda (\G )$--rank of $\CC _p (E/\F )$ is isogeny invariant.
\end{Remark}
\begin{proof}
This follows immediately from Theorem \ref{aa}. Since we are now assuming the $\CC _p (E_i /\F )$ are $\Lambda (\G )$--torsion (and thus also $\Lambda (G_L )$--torsion for all finite extensions, $L$ of $F$) we know from Proposition \ref{bb} that
\begin{equation}\label{406}
\begin{array}{c}
\mu _L ({\widehat{\cok ( \phi _1 )}}) - \mu _L (\CC _p (E_2 /\F)) 
+ \mu _L (\CC _p (E_1 /\F)) - \mu _L ({\widehat{\kr ( \phi _1 )}}) = 0.
\end{array}
\end{equation}
Replacing $F$ by $L$ and $\G $ by $G_L $ in Theorem \ref{aa}, and substituting the formula thus resulting from (\ref{212}) into (\ref{406}), gives the Corollary.{\\ \nopagebreak
\hspace*{\fill}$\Box$}
\end{proof}

\begin{ex}
We consider the isogeny class of non-CM elliptic curves of conductor 
11 containing the curve $X_0(11)$ considered above, (\ref{70}). Take 
$E_1 $ to be the curve $ X_1(11)$, given by a minimal Weierstra\ss $ \; $
equation 
$y^2 + y = x^3 - x^2$, and let $E_2 $ be $ X_0 (11)$. 
 They are related via the isogeny
$\phi : E_1 \rightarrow E_2 $ with kernel $A \cong \z /5 $. As remarked in the example at the end of \S 2, it 
follows from Lang and Trotter's determination of $\G $ in this 
case that $\Gal (\q (E_{i,5^{\infty }}) / \q (\mu _5 ) )$ is pro-5, 
\cite{LaTr}. Letting $L$ denote any finite extension of 
$\q (\mu _5 )$ contained in $\q (E_{i,5^{\infty }})$ and so 
$G_L $ is the group $\Gal (E_{i,5^{\infty }}/L)$, 
it is shown in \cite{CoHo99}, 
\S 6 and 7, that both $\CC _5 (E_i /\F )$ are $\Lambda (G_L )$--torsion, 
and thus by Proposition 3.4 of \cite{CoHo99} Assumption {\bf{L}} holds. 
Then Corollary \ref{tt} applies.
Since $L$ contains $\q (\mu _5 )$, all Archimedean places of $L$ are complex. Then 
\begin{equation}\label{403}
\begin{array}{c}
\underset{\omega \mid \infty }{\sum} \, \ord _5 \big( \# A (L_{\omega }) \big) = \frac{1}{2} \mid L:\q \mid ,\\
\mid L:\q \mid \ord _5 (\# A) = \mid L:\q \mid .
\end{array}
\end{equation}
Turning now to local considerations, the kernel of the isogeny $\phi : E_1 \longrightarrow E_2 $ injects into $\tilde{E}_{1,5^{\infty }}$ and so $\tilde{A}_{\omega } \cong \z /5 $ also, for each prime $\omega $ of $L$ above 5. Thus
\begin{equation}\label{404}
\underset{\omega \mid 5 }{\sum }\, \ord _5 \big( \mid \# \tilde{A}_{\omega } \mid _{\omega } \big) = - \mid L:\q \mid .
\end{equation}
It is also explained in \cite{CoHo99}, Corollary 7.3, that $\CC _5 (E_1 /\F )$ is a finitely generated $\Lambda (H)$--module for $H = \Gal (\F /L^{cyc})$. Thus, by Lemma \ref{yy}, $\mu _L (\CC _5 (E_1 /\F )) = 0.$ We are using here again the notation $\mu _L $ defined in Corollary \ref{tt}.
From this, the calculations (\ref{403}) and (\ref{404}) above and Corollary \ref{tt} we see that as a $\Lambda (G_L )$--module
\begin{equation}
\mu _L (\CC _5 (E _2 /\q (E_{2,5^{\infty }}))) = {{\textstyle{\frac{1}{2}}}} \mid L:\q \mid .
\end{equation}
\end{ex}
Similar calculations can be made also for the third curve in the isogeny class.



\end{document}